\documentclass[12pt,english]{article}
\usepackage{amssymb}
\usepackage{amsfonts}
\usepackage[english]{babel}

\newcommand{\Z}{{\mathbb Z}}
\newcommand{\N}{{\mathbb N}}
\newcommand{\R}{{\mathbb R}}
\newcommand{\av}[1]{\langle #1 \rangle}

\newtheorem{theorem}{Theorem}[section]
\newtheorem{lemma}[theorem]{Lemma}
\newtheorem{e-proposition}[theorem]{Proposition}

\newtheorem{e-definition}[theorem]{Definition\rm}

\begin{document}

\begin{titlepage}
\vskip 0.5cm
\vskip 1.5cm
\begin{center}
{\Large{\bf
Ces\`aro asymptotics for the orders of $SL_k(\mathbb{Z}_n)$ and $GL_k(\mathbb{Z}_n)$ as $n\to\infty$
}}
\end{center}
\vskip 0.8cm
\centerline{Alexey G. Gorinov\footnote{e-mail: {\tt gorinov@math.jussieu.fr}}
and Sergey V. Shadchin\footnote{e-mail: {\tt chadtchi@ihes.fr}}}
\vskip 0.9cm
\centerline{${}^{1}$\sl\small Universit\'e Paris 7, U.F.R. de Math\'ematiques, 2 place Jussieu, F-75251, France \,}
\centerline{${}^{2}$\sl\small IHES, Bures-sur-Yvette, route de Chartres, F-91140, France \,}

\vskip 1.25cm

\begin{abstract}
Given an integer $k>0$, our main result states that the sequence of 
orders of the groups $SL_k(\Z_n)$ (respectively, of the groups $GL_k(\Z_n)$)
is Ces\`aro equivalent as $n\to\infty$ to the sequence 
$C_1(k) n^{k^2-1}$ (respectively,
$C_2(k)n^{k^2}$), where the coefficients $C_1(k)$ and $C_2(k)$ 
depend only on $k$; we give explicit formulas for $C_1(k)$ and $C_2(k)$.
This result generalizes the theorem
(which was first published by I. Schoenberg) that says
that the Euler function $\varphi (n)$ is Ces\`aro equivalent to 
$\frac{6}{\pi^2} n$.
We present some experimental facts related to
the main result.
\end{abstract}

\end{titlepage}


\section{Introduction}
The article is organized as follows: in section \ref{sec:theorem} we introduce some notation and formulate our main result. Then, in section \ref{sec:proof}, we prove this result. Finally, in section \ref{sec:experimental} we discuss some interesting related facts, in particular, some experimental facts.

The authors wish to thank V. I. Arnold for proposing this interesting problem and for useful discussions.

\section{The main theorem}
\label{sec:theorem}
Two sequences of real numbers $(x_n)_{n\in\N}$ and $(y_n)_{n\in\N}$ are said to be
{\it Ces\`aro equivavent}, if $$\lim_{n\to\infty}\frac{x_1+\cdots +x_n}
{y_1+\cdots +y_n}=1.$$ For any finite set $X$ we shall denote by
$\# (X)$ the cardinality of $X$. We shall use the symbol
$\prod_p$ to denote the product over all prime numbers.

Our main result is the following theorem.
\begin{theorem}
\label{main}
For any fixed integer $k>0$ the sequence $\left(\# (SL_k(\Z_n))\right)_{n\in\N}$ $($respectively, the sequence $\left(\# (GL_k(\Z_n))\right)_{n\in\N}$$)$ is Ces\`aro
equivalent as $n\to\infty$ to $C_1(k)n^{k^2-1}$ $($respectively, $C_2(k)n^{k^2}$$)$, where ${C_1(1)=1}, C_2(1)=\prod_p\left(1-\frac{1}{p^2}\right)$, and
for any $k>1$ we have
$$C_1(k)=\prod_p\left( 1-\frac{1}{p}\left(1-
\prod_{i=2}^k \left(1-\frac{1}{p^i}\right)\right)\right),$$
$$ C_2(k)=\prod_p\left(
1-\frac{1}{p}\left(1-
\prod_{i=1}^k\left(1-\frac{1}{p^i}\right)\right)\right).$$
\end{theorem}

{\bf Remark.} In particular, $\# (GL_1(\Z_n))$ and $\# (SL_2(\Z_n))$
are Ces\`aro equivalent to $\frac{n}{\zeta(2)}$ and
$\frac{n^3}{\zeta(3)}$ respectively. We do not know if 
the asymptotics given by Theorem \ref{main} can be expressed
in terms of values of the Riemann zeta-function
(or any other remarkable function)
at algebraic points in any of the other cases. 

To the best of our knowledge, the fact that the Euler 
function $\varphi(n)=\#( GL_1(\Z_n) )$ is Ces\`aro equivalent to
$n\frac{6}{\pi^2}$ was first
published in \cite{schonberg} by I. Schoenberg, who attributes the 
result to J. Schur. This result was probably
already known to Gauss. An explicit formula for the cumulative 
distribution function of the sequence 
$\left(\frac{\varphi(n)}{n}\right)_{n\in\N}$ is 
given in \cite{venkov} by B. A. Venkov.

\section{Proof of Theorem \ref{main}}
\label{sec:proof}
Let us first recall the explicit formulas 
for $\# (SL_k(\Z_n))$ and $\# (GL_k(\Z_n))$.
For any positive integer $k$ denote by
$\tilde\varphi_k$\footnote{This notation can be explained as 
follows: the function $\tilde\varphi_k$
generalizes the function
$n\mapsto\varphi(n)/n=\tilde\varphi_1(n)$.} the map $\N\to\R$ given
by the formula 
$$\tilde\varphi_k(p_1^{l_1}\cdots p_m^{l_m})=\left(1-\frac{1}{p_1^k}
\right)
\cdots \left(1-\frac{1}{p_m^k}\right)$$ 
(here $p_1,\ldots ,p_m$ are pairwise distinct primes).

\begin{lemma}
\label{expl} We have $$\# (GL_1(\Z_n))=n\tilde\varphi_1(n),$$ and for any integer $k>1$ we have
$$\# (SL_k(\Z_n))=n^{k^2-1}\tilde\varphi_2(n)\cdots\tilde\varphi_k(n),$$
$$\# (GL_k(\Z_n))=n^{k^2}\tilde\varphi_1(n)\cdots\tilde\varphi_k(n).$$
\end{lemma}

{\it Proof of Lemma \ref{expl}.} If $p$ is prime, then $$\# (GL_k(\Z_p))=
(p^k-1)(p^k-p)\cdots (p^k-p^{k-1}).$$

The formula for $\# (GL_k(\Z_n))$ follows
from the existence of the homomorphisms $$GL_k(\Z_{p^l})\to GL_k(\Z_{p^{l-1}})$$
(for any prime $p$ and for any positive integer $l$),
and from the fact that $$GL_k(\Z_{nm})=GL_k(\Z_{n})\times
GL_k(\Z_{m}),$$ if $m$ and $n$ are positive coprime integers.
The formula for $\# (SL_k(\Z_n))$ is obtained from the formula for $\# (GL_k(\Z_n))$ using
the determinant homomorphism.$\diamondsuit$

Now let us calculate the limits of the averages of the sequences 
$(\tilde\varphi_1(n)\cdots\tilde\varphi_k(n))_{n\in\N}$
and $(\tilde\varphi_2(n)\cdots\tilde\varphi_k(n))_{n\in\N}$. More generally, 
let $\ell$ be a finite (nonempty)
ordered collection of positive integers: 
$\ell=(i_1,\ldots,i_l)$. For any $n\in\N$ set $$\tilde\varphi_\ell(n)=\tilde\varphi_{i_1}(n)
\cdots\tilde\varphi_{i_l}(n).$$ For any sequence $x=(x_n)_{n\in\N}$ denote by $\av{x}$ the Ces\`aro 
limit of $x$, i.e.,  $$\lim_{n\to\infty}\frac{1}{n}\sum_{m=1}^nx_m = \av{x}.$$

\begin{theorem}\label{formula}

For any $\ell=(i_1,\ldots,i_l)$ the limit $\av{\tilde\varphi_\ell}$
exists and is equal to $\prod_p f_\ell\left(\frac{1}{p}\right)$, 
where $$f_\ell(t)=1-t(1-\prod_{j=1}^l (1-t^{i_j})).$$
\end{theorem}

{\it Sketch of a proof of Theorem \ref{formula}.}
We shall first give an informal proof of the theorem; we shall then
show what changes should be made to make our informal proof rigorous.

The idea of the proof of Theorem \ref{formula} is to give a 
probabilistic interpretation to some complicated expressions (such as
$\frac{1}{n}\sum_{m=1}^n \tilde\varphi_\ell(m)$). This idea goes back to Euler.

Let us note that for any positive integer $q$ the ``probability''
that a ``random'' positive integer is a not a multiple of $q$ is $$1-\frac{1}{q}.$$ 
If $q_1$ and $q_2$
are coprime integers, the events ``$r$ is not divisible by $q_1$'' and 
``$r$ is not divisible by $q_2$'' ($r$ being a ``random'' positive
integer) are independent, which implies that
for any positive integers $m,k$ the expression 
$\tilde\varphi_k(m)$ is the ``probability'' that a
``randomly chosen'' positive integer is not divisible
by $k$-th powers of the prime divisors of $m$.

Analogously, for any fixed positive integer $m$ the expression $\tilde\varphi_\ell(m)$ can be seen as the 
``probability'' to find an element $$(x_1,\ldots ,x_l)\in\N^l$$ that satisfies the following conditions:
$x_1$ is not divisible by the $i_1$-th powers of the 
prime factors of $m$, $x_2$ is not divisible by
the $i_2$-th powers of the prime factors of $m$ etc.

Using the total probability formula, we obtain that
$\frac{1}{n}\sum_{m=1}^n \tilde\varphi_\ell(m)$ is the ``probability''
that a ``random'' element of the set $$\{(x_0,x_1\ldots,x_l)|x_0,\ldots,x_l\in\N, x_0\leq n\}$$ satisfies the following condition:
any $x_j, j=1,\ldots,l$ is not divisible by the
$i_j$-th powers of the prime divisors of $x_0$.
So the limit $\av{\tilde\varphi_\ell}$ 
is the ``probability'' of the limit event, which can be
described
as the intersection for all prime $p$ of the following events:
``either ($x_0$ is not divisible by $p$), or (none of $x_j, j=1,\ldots ,l$ is divisible by $p^{i_j}$)''.
These events are independent, and the ``probability'' of each of 
them is $$f_\ell\left(\frac{1}{p}\right)=
1-\frac{1}{p}\left(1-\prod_{j=1}^l\left(1-\frac{1}{p^{i_j}}\right)\right).$$
This gives the desired expression for $\av{\tilde\varphi_\ell}$.

This idea is formalized as follows. Let $l$ be a 
positive integer, and let $A$ and $B$ be subsets of $\N^l$ such that
there exists $$\lim_{k\to\infty}\frac{\# (A\cap B\cap C_k)}{\# (B\cap C_k)},$$ 
where $C_k=\{(x_1,\ldots,x_l)\in\N^l|x_1\leq k,\ldots,x_l\leq k\}$. 
This limit will be called
the {\it density} of $A$ in $B$ and will be denoted by $p_B(A)$. 
For any $B\subset\N^l$ the correspondence 
$B\supset A\mapsto p_B(A)$ defines a measure on
$B$\footnote{Unfortunately, this measure is not $\sigma$-additive, 
which is why we prefer to speak rather of densities than of probabilities.}.

Using the same argument as above (and replacing 
``probabilities'' with ``densities'' and ``events'' with ``sets''), 
we can represent $\frac{1}{n}\sum_{m=1}^n\tilde\varphi_\ell(m)$ as the density of a certain subset of the
set $$\{(x_0,x_1\ldots,x_l)|x_0,\ldots,x_l\in\N, x_0\leq n\}$$. 
This interpretation does not allow us to pass immediately to the 
limit as $n\to\infty$, but
it enables us to write the following combinatorial 
formula for $\frac{1}{n}\sum_{m=1}^n\tilde\varphi_\ell(m)$. Define the sequence $(a_k)_{k\in\N}$ by the formula
$$\sum_{k=1}^\infty a_kt^k=1-\prod_j(1-t^{i_j}).$$ We have 
$$\frac{1}{n}\sum_{m=1}^n \tilde\varphi_\ell(m)=1+
\sum_{r=2}^\infty \frac{1}{r}(-1)^{pr(r)}a(r) b_{r,n},$$
where for any $r=p_1^{\alpha_1}\cdots p_s^{\alpha_s}$ we define
$$pr(r)=s,\;\;\; a(r)=a_{\alpha_1}\cdots a_{\alpha_s}, \;\;\; 
b_{r,n}=\left[\frac{n}{p_1\cdots p_s}\right]
\frac{1}{n}$$ (in particular, $a(r)=0$, 
if $\max\{\alpha_1,\ldots ,\alpha_s\}>i_1+\cdots +i_l$). 
Now let us note that this expression
has the form $\sum_{k=1}^\infty b'_{k,n}c_r$, where 
$c_k$ is the $k$-th term of the absolutely
convergent series obtained by multiplying out the product 
$$\prod_p \left(1-\frac{1}{p}
\left(1-\prod_{j=1}^l\left(1-\frac{1}{p^{i_j}}\right)\right)\right),$$
and every $b'_{k,n}$ has the form
$$\frac{p_1\cdots p_s}{n}\left[\frac{n}{p_1\cdots p_s}\right].$$ 
We have $0\leq b'_{k,n}\leq 1$ for any $k,n$, and the limit $\lim_{n\to\infty} b'_{k,n}$ is equal to 1 for any $k$. This implies 
Theorem \ref{formula}.
$\diamondsuit$

Theorem \ref{main} can be obtained from Theorem \ref{formula}, 
from Lemma \ref{expl} and from the following lemma.

\begin{lemma}\label{stepen}
Let $(x_n)_{n\in\N}$ be a sequence of real numbers, and suppose 
that $\av{x}$ exists. Then, for any 
nonnegative integer $k$, we have 
$$\lim_{n\to\infty}\frac{x_1+2^kx_2+\cdots+ n^kx_n}{1+2^k+\cdots+n^k}=\av{x}.$$
\end{lemma}

{\it Proof of Lemma \ref{stepen}.} The proof is by induction on $k$. 
If $k=0$, there is nothing to prove. Suppose Lemma \ref{stepen} holds for
some $k$. For any sequence $y=(y_n)_{n\in\N}$ set 
$$S_n^k[y]=y_1+2^ky_2+\cdots+n^ky_n.$$ We have 
$$S_n^k[x]=n^{k+1}\left(\frac{\av{x}}{k+1}+\varepsilon_n,\right)$$ where 
$(\varepsilon_n)_{n\in\N}$ is a sequence such that 
$\lim_{n\to\infty}\varepsilon_n=0$.
Note that for any sequence $y=(y_n)_{n\in\N}$ we have $$S_n^{k+1}[y] = n S_n^k[y] - \sum_{m=1}^{n-1}S_m^k[y].\leqno{(*)}$$
Thus, we can write 
$$S_n^{k+1}[x] =
\frac{\av{x}}{k+1}\left( n^{k+2} - \sum_{m=1}^{n-1}m^{k+1}\right) 
+ \varepsilon_n n^{k+2} - S_{n-1}^{k+1}[\varepsilon]$$. 

We have 
$$\lim_{n\to\infty}\frac{S_{n-1}^{k+1}[\varepsilon]}{n^{k+2}} = 0,$$ and hence 
$$\lim_{n\to\infty}\frac{S_n^{k+1}[x]}{n^{k+2}} =
\frac{\av{x}}{k+1}\left(1-\frac{1}{k+2}\right) = \frac{\av{x}}{k+2},$$
which implies the statement of Lemma \ref{stepen}. $\diamondsuit$

\section{Convergence rates and the distribution of the 
values of $\tilde\varphi_\ell$}
\label{sec:experimental}
Let $\ell$ be a finite (nonempty) ordered collection of positive
integers: $\ell=(i_1,\ldots, i_l)$.
In this section we briefly discuss the convergence rate of the sequences
$$\left(\frac{1}{n^{s+1}}\sum_{k=1}^n k^s\tilde{\varphi}_\ell(k)\right)_{n\in\N}$$ for different fixed $s\in\N$ and the 
distribution of the values of the function
$\tilde\varphi_\ell$.

Set $$\Phi_{\ell} = \lim_{n\to\infty} \frac{1}{n}
\sum_{k=1}^n\tilde{\varphi}_{\ell}(k) = \prod_p f_{\ell}(\frac{1}{p}),$$ 
$$\xi_{\ell,s}(n) =\frac{1}{n^s} 
\left(\sum_{k=1}^n k^s \tilde{\varphi}_{\ell}(n) - 
\frac{n^{s+1}}{s+1}\Phi_{\ell}\right).$$
It follows immediately from these definitions that 
$$\sum_{k=1}^n k^s \tilde{\varphi}_{\ell}(k) = 
\frac{n^{s+1}}{s+1}\Phi_{\ell} + n^s \xi_{\ell,s}(n).$$

We have carried out computer experiments for $\ell=(1)$, $(2)$ 
and $(1,2)$ and for $s=0,...,3$. In all these cases the 
sequence $\xi_{\ell,s}=(\xi_{\ell,s}(n))_{n\in\N}$ seems to be bounded.
This means that the first relative correction to $\Phi_{\ell}$ should be of
order $O\left(\frac{1}{n}\right)$. Computer tests suggest that the 
sequence $\xi_{\ell,s}$ has no limit, but has a Ces\`aro limit, 
which is equal to $\frac{1}{2}\Phi_{\ell}$ for every nonnegative integer $s$. At this moment we are able to prove only the following weaker statement.

\begin{theorem}\label{the:averages}
If $\av{\xi_{\ell,0}}$ exists, then for all integers $s>0$ the limit
$\av{\xi_{\ell,s}}$ exists and is equal to $\frac{1}{2}\Phi_{\ell}$.
\end{theorem}

{\it Proof of Theorem \ref{the:averages}.} Set 
$$\eta_{\ell,s}(n)=
\frac{1}{n^s}\sum_{k=1}^nk^s\left(\tilde\varphi_\ell (k)-\Phi_\ell\right).$$ 
Note that $\xi_{\ell,0}=\eta_{\ell,0}$, hence
$\av{\eta_{\ell,0}}$ exists. Using formula $(*)$ we get $$\eta_{\ell,s+1}(n)=\eta_{\ell,s}(n)-\frac{1}{n^{s+1}}\sum_{k=1}^n k^s\eta_{\ell,s}(k).$$
Hence we obtain using Lemma~\ref{stepen} that
$\av{\eta_{\ell,s+1}}=\frac{s}{s+1}\av{\eta_{\ell,s}}$ for any integer $s\geq 0$. Thus, $\av{\eta_{\ell,s}}=0$ for any integer $s>0$.

For any integer $s\geq 1$ we have 
$$\sum_{k=1}^n k^s = \frac{n^{s+1}}{s+1} +
\frac{1}{2}n^s + O\left(n^{s-1}\right).$$ 
Hence we get the following 
relation: 
$$\xi_{\ell,s}(n)~=~\eta_{\ell,s}(n)~+~\frac{1}{2}\Phi_{\ell}+O\left(\frac{1}{n}\right),$$ 
which implies that
$\av{\xi_{\ell,s}} = \frac{1}{2}\Phi_{\ell}$.
The theorem is proven. $\diamondsuit$

Let us now consider the distribution of the values of the 
function $\tilde{\varphi}_{\ell}$.
Using the argument from \cite[\S 5]{schonberg}, one can prove that for any $t\in [0,1 ]$ the limit
$$\lim_{n\to\infty}\frac{1}{n}\#\{k\in\N|k\leq n,\tilde\varphi_\ell (k)\leq t\}$$ exists, and that the function $F_\ell$ defined by the formula $$F_\ell (t)=\lim_{n\to\infty}\frac{1}{n}\#\{k\in\N|k\leq n,\tilde\varphi_\ell (k)\leq t\}$$
is continuous (I.~Schoenberg considers only the case $\ell=(1)$, but
his argument can be easily extended to the case of
an arbitrary $\ell$). The function $F_\ell$ is the analogue of the cumulative distribution function in probability theory.
Given a nonnegative integer $s$, the $s$-th moment of $F_\ell$ is defined as follows: 
$$\mu_{\ell,s} = \int_0^1t^s dF_\ell(t).$$
It is easy to prove (see \cite[Satz I]{schonberg}) that $\mu_{\ell,s}=\av{{(\tilde{\varphi}_{\ell})}^s}$.
Due to Theorem \ref{formula}, we have $\mu_{\ell,s} = \Phi_{\ell^s}$ where $\ell^s$ is the
following collection of positive integers: 
$\ell^s=( i_1,i_1,...,i_1 \mbox{ ($s$ times) }, i_2,i_2,....,i_2 
\mbox{ ($s$ times) },...)$. 

The Fourier series for $F_{\ell}(t)$ is equal to $\sum_{n\in\Z}u_ne^{2\pi
int}$, where $$u_0=1-\Phi_\ell=\frac{1}{2}-\sum_{k\neq 0}\frac{\av{e^{-2\pi ik\tilde\varphi_\ell}}}{2\pi ik}$$ (the sum of the series in the latter formula is to be taken in Ces\`aro sence),
and the Fourier coefficients $u_n$ for
$n\neq 0$ can be calculated using either the formula
$$u_n=-\sum_{m=1}^{\infty}\frac{(-2\pi i n)^{m-1}}{m!}\Phi_{\ell^m},$$ or the formula $$u_n=\frac{1}{2\pi in}(\av{e^{-2\pi in\tilde\varphi_\ell}}-1).$$

Since $F_\ell $ is continuous, its Fourier series converges in Ces\`aro sence to $F_\ell$ uniformly on every compact subset of the open interval $(0,1)$.


\begin{thebibliography}{99}
\bibitem{schonberg} I. Schoenberg, \"Uber die asymptotische 
Verteilung reeller Zahlen mod 1, Math. Z. 28 (1928), 171-199.
\bibitem{venkov} B. A. Venkov, On a certain monotonic function, 
Uchenye Zapiski Leningrad State Univ., Math. Ser. 111, fasc. 16 (1949), 
3-19 (in Russian). 

\end{thebibliography}
\end{document}